\documentclass[a4paper]{amsart}

\newtheorem{theorem}{Theorem}[section]
\newtheorem{corollary}[theorem]{Corollary}

\newtheorem{lemma}[theorem]{Lemma}
\theoremstyle{remark}
\newtheorem{remark}{Remark}[section]

\theoremstyle{definition}

\theoremstyle{definition}

\begin{document}

\title[Combined dynamic Gr\"{u}ss inequalities on time scales]{%
Combined dynamic Gr\"{u}ss inequalities\\ on time scales}



\author[M. R. Sidi Ammi\ and\ D. F. M. Torres]{Moulay Rchid Sidi Ammi \and Delfim F. M. Torres}

\address{Department of Mathematics,
University of Aveiro,
3810-193 Aveiro, Portugal}

\email{sidiammi@ua.pt}
\email{delfim@ua.pt}

\thanks{Research partially supported by the
{\it Centre for Research on Optimization and Control} (CEOC)
from the {\it Portuguese Foundation for Science and Technology} (FCT),
cofinanced by the European Community Fund FEDER/POCI 2010; and
FCT postdoc fellowship SFRH/BPD/20934/2004.}

\maketitle

%---------------------------------------------------

\begin{abstract}
We prove a more general version of the Gr\"{u}ss inequality
by using the recent theory of combined dynamic derivatives
on time scales and the more general notions
of diamond-$\alpha$ derivative and integral.
For the particular case when $\alpha = 1$,
one gets a delta-integral Gr\"uss inequality on time scales;
for $\alpha = 0$ a nabla-integral Gr\"uss inequality.
If we further restrict ourselves by fixing the time
scale to the real (or integer) numbers, then the 
standard continuous (discrete) inequalities are obtained.
\end{abstract}

%---------------------------------------------------

\section{Introduction}

The Gr\"{u}ss inequality is of great interest in differential
and difference equations as well as many other areas of
mathematics \cite{d:Cerone,d:dragomir,d:Liu,mit,d:mit}.
The classical inequality was proved by G.~Gr\"{u}ss in 1935 
\cite{gruss}: if $f$ and $g$ are two continuous functions 
on $[a,b]$ satisfying $\varphi \le f(t) \le \Phi$ 
and $\gamma \le g(t) \le \Gamma$ for all $t \in [a,b]$, then
\begin{equation}
\label{eq:class:Gruss}
\left| \frac{1}{b-a}\int_{a}^{b}  f(t) g(t) dt
- \frac{1}{(b-a)^{2}}\int_{a}^{b} f(t) dt  \int_{a}^{b} g(t) dt
\right| \leq \frac{1}{4} (\Phi-\varphi)(\Gamma-\gamma) \, .
\end{equation}
The literature on Gr\"uss type inequalities is now vast,
and many extensions of the classical inequality \eqref{eq:class:Gruss}
have been intensively studied by many authors in the XXI century
\cite{d:abra,bohner,dragomir,d:imp,d:Pachpatte,d:gti,d:gti:d,d:la}.
Here we are interested in Gr\"uss type inequalities on time scales.

The analysis on time scales is a relatively new area of mathematics
that unify and generalize discrete and continuous theories.
Moreover, it is a crucial tool in many computational and numerical
applications. The subject was initiated by Stefan Hilger \cite{h2,h3}
and is being applied to many different fields
in which dynamic processes can be described
with discrete or continuous models:
see \cite{abra,comZbig,b1,b2,with-Rui} and references therein.
One of the important subjects being developed within
the theory of time scales is the study
of inequalities \cite{abp,srd,SteffensenIneq,R:T:JIPAM}.

Recently, Bohner and Matthews \cite{bohner} proved a
time scales version of the Gr\"{u}ss inequality
\eqref{eq:class:Gruss} by using the delta integral. Roughly
speaking, the main result in \cite{bohner} asserts
that \eqref{eq:class:Gruss} continues to be valid 
on a general time scale by substituting the Riemann integral 
by the  $\Delta$-integral. The main objective of this paper 
is to use the more general diamond$-\alpha$ dynamic integral.

In 2006, a combined dynamic derivative $\diamondsuit_{\alpha}$ was
introduced as a linear combination of $\Delta$ and $\nabla$ dynamic
derivatives on time scales \cite{sfhd}.
The diamond--$\alpha$ derivative reduces
to the standard $\Delta$ derivative for $\alpha =1$ and to the standard
$\nabla$ derivative for $\alpha =0$. On the other hand, it
represents a weighted dynamic derivative formula on any uniformly
discrete time scale when $\alpha =\frac{1}{2}$. We refer the reader
to \cite{DorotaControlo2008,Rogers,Sheng,sfhd} for an account 
of the calculus with diamond-$\alpha$ dynamic derivatives.

In Section~\ref{sec:Prel} we briefly review the
necessary definitions and calculus on time scales.
Our results are given in Section~\ref{sec:MR}.

%---------------------------------------------------

\section{Preliminaries}
\label{sec:Prel}

A time scale $\mathbb{T}$ is an arbitrary nonempty closed subset of
the real numbers. The time scales calculus was initiated 
by S.~Hilger in his PhD thesis in order to unify discrete and
continuous analysis \cite{h2,h3}. Let $\mathbb{T}$ 
be a time scale with the topology that 
it inherits from the real numbers. For $t \in \mathbb{T}$, 
we define the forward jump operator 
$\sigma: \mathbb{T} \rightarrow \mathbb{T}$ by 
$$\sigma(t)= \inf \left \{s \in \mathbb{T}: s >t \right\} \, ,$$ 
and the backward jump operator  $\rho: \mathbb{T} \rightarrow \mathbb{T}$ 
by $$\rho(t)= \sup \left \{s \in \mathbb{T}: s < t \right \} \, .$$

If $\sigma(t) > t$ we say that $t$ is right-scattered, 
while if $\rho(t) < t$ we say that $t$ is left-scattered. 
Points that are simultaneously right-scattered and 
left-scattered are said to be isolated. 
If $\sigma(t)=t$, then $t$ is called right-dense; 
if $\rho(t)=t$, then $t$ is called left-dense. 
Points that are right-dense and left-dense 
at the same time are called dense. 
The mappings $\mu, \nu: \mathbb{T}
\rightarrow [0, +\infty)$ defined by $\mu(t):=\sigma(t)-t$ 
and $\nu(t):=t-\rho(t)$ are called, respectively, 
the forward and backward graininess function.

Given a time scale $\mathbb{T}$,
we introduce the sets $\mathbb{T}^{k}$, $\mathbb{T}_{k}$, and
$\mathbb{T}^{k}_{k}$ as follows. If $\mathbb{T}$ has a left-scattered
maximum $t_{1}$, then $\mathbb{T}^{k}= \mathbb{T}-\{t_{1} \}$,
otherwise $\mathbb{T}^{k}= \mathbb{T}$. If  $\mathbb{T}$ has a
right-scattered minimum $t_{2}$, then $\mathbb{T}_{k}=
\mathbb{T}-\{t_{2} \}$, otherwise $\mathbb{T}_{k}= \mathbb{T}$.
Finally, $\mathbb{T}_{k}^{k}= \mathbb{T}^{k} \bigcap
\mathbb{T}_{k}$.

Let $f: \mathbb{T}\rightarrow \mathbb{R}$ be a real 
valued function on a time scale $\mathbb{T}$. 
Then, for $t \in \mathbb{T}^{k}$, we define $f^{\Delta}(t)$
to be the number, if one exists, such that for all $\epsilon >0$,
there is a neighborhood $U$ of $t$ such that for all $s \in U$,
$$
\left|f(\sigma(t))- f(s)-f^{\Delta}(t)(\sigma(t)-s)\right| \leq \epsilon
|\sigma(t)-s|.
$$
We say that $f$ is delta differentiable on $\mathbb{T}^{k}$
provided $f^{\Delta}(t)$ exists for all $t \in \mathbb{T}^{k}$.
Similarly, for $t \in \mathbb{T}_{k}$ we define $f^{\nabla}(t)$ to
be the number, if one exists, such that for all $\epsilon >0$,
there is a neighborhood $V$ of $t$ such that for all $s \in V$
$$
\left|f(\rho(t))- f(s)-f^{\nabla}(t)(\rho(t)-s)\right|
\leq \epsilon |\rho(t)-s|.
$$
We say that $f$ is nabla differentiable on $\mathbb{T}_{k}$,
provided that $f^{\nabla}(t)$ exists for all $t \in \mathbb{T}_{k}$.

For $f:\mathbb{T}\rightarrow \mathbb{R}$ we
define the function $f^{\sigma}: \mathbb{T}\rightarrow \mathbb{R}$
by $f^{\sigma}(t)=f(\sigma(t))$ for all $t \in \mathbb{T}$, that is,
$f^{\sigma}= f\circ \sigma$. Similarly, we
define the function $f^{\rho}: \mathbb{T}\rightarrow \mathbb{R}$ by
$f^{\rho}(t)=f(\rho(t))$ for all $t \in \mathbb{T}$, that is,
$f^{\rho}= f\circ \rho$.

A function $f: \mathbb{T} \rightarrow \mathbb{R} $ is called
rd-continuous, provided it is continuous at all right-dense points
in $\mathbb{T}$ and its left-sided limits finite at all left-dense
points in $\mathbb{T}$. A function $f: \mathbb{T} \rightarrow \mathbb{R} $ is called
ld-continuous, provided it is continuous at all left-dense points in
$\mathbb{T}$ and its right-sided limits finite at all right-dense
points in $\mathbb{T}$.

A function $F: \mathbb{T} \rightarrow \mathbb{R} $ is called a delta
antiderivative of $f: \mathbb{T} \rightarrow \mathbb{R}$, provided
$F^{\Delta}(t)=f(t)$ holds for all $t \in \mathbb{T}^{k}$. Then the
delta integral of $f$ is defined by $$\int^b_a f(t)\Delta
t=F(b)-F(a) \, .$$

A function $G: \mathbb{T} \rightarrow \mathbb{R} $ 
is called a nabla antiderivative of $g: \mathbb{T} 
\rightarrow \mathbb{R}$, provided $G^{\nabla}(t)=g(t)$ 
holds for all $t \in \mathbb{T}_{k}$. Then the
nabla integral of $g$ is defined by $\int^b_a g(t)\nabla
t=G(b)-G(a)$.
For more details on time scales one can see \cite{abra,b1,b2}.

Now we introduce the diamond-$\alpha$ derivative
and integral, referring the reader to 
\cite{DorotaControlo2008,Rogers,Sheng,sfhd} 
for more on the associated calculus.

Let $\mathbb{T}$ be a time scale and $f$ differentiable on
$\mathbb{T}$ in the $\Delta$ and $\nabla$ senses. For $t \in
\mathbb{T}$, we define the diamond-$\alpha$ dynamic derivative
$f^{\diamondsuit_{\alpha}}(t)$ by
$$
f^{\diamondsuit_{\alpha}}(t)= \alpha
f^{\Delta}(t)+(1-\alpha)f^{\nabla}(t), \quad 0 \leq \alpha \leq 1.
$$
Thus, $f$ is diamond-$\alpha$ differentiable if and only if $f$ 
is $\Delta$ and $\nabla$ differentiable. The diamond-$\alpha$ 
derivative reduces to the standard  $\Delta$ derivative for $\alpha =1$, 
or the standard $\nabla$ derivative for $\alpha =0$. 
On the other hand, it represents a ``weighted derivative'' 
for $\alpha \in (0,1)$. Diamond-$\alpha$ derivatives have shown 
in computational experiments to provide efficient 
and balanced approximation formulas, leading to 
the design of more reliable numerical methods \cite{Sheng,sfhd}.

Let $f, g: \mathbb{T} \rightarrow \mathbb{R}$ be diamond-$\alpha$
differentiable at $t \in \mathbb{T}$. Then,

\begin{itemize}

\item[(i)] $f+g: \mathbb{T} \rightarrow \mathbb{R}$ is diamond-$\alpha$
differentiable at $t \in \mathbb{T}$ with
$$ (f+g)^{\diamondsuit^{\alpha}}(t)=
(f)^{\diamondsuit^{\alpha}}(t)+(g)^{\diamondsuit^{\alpha}}(t).
$$

\item[(ii)] For any constant $c$, $cf:  \mathbb{T} \rightarrow \mathbb{R}$
 is diamond-$\alpha$ differentiable at $t \in \mathbb{T}$ with
$$
(cf)^{\diamondsuit^{\alpha}}(t)= c(f)^{\diamondsuit^{\alpha}}(t).
$$

\item[(ii)]  $fg: \mathbb{T} \rightarrow \mathbb{R}$ 
is diamond-$\alpha$ differentiable at $t \in \mathbb{T}$ with
$$
(fg)^{\diamondsuit^{\alpha}}(t)= (f)^{\diamondsuit^{\alpha}}(t)g(t)+
\alpha f^{\sigma}(t)(g)^{\Delta}(t) +(1-\alpha)
f^{\rho}(t)(g)^{\nabla}(t).
$$

\end{itemize}

Let $a, t \in  \mathbb{T}$, and $h: \mathbb{T} \rightarrow
\mathbb{R}$. Then, the diamond-$\alpha$ integral from $a$ to $t$ of $h$ is defined by
$$
\int_{a}^{t}h(\tau) \diamondsuit_{\alpha} \tau = \alpha
\int_{a}^{t}h(\tau) \Delta \tau +(1- \alpha) \int_{a}^{t}h(\tau)
\nabla \tau, \quad 0 \leq \alpha \leq 1.
$$
We may notice the absence of an anti-derivative
for the $ \diamondsuit_{\alpha}$ combined derivative.
For $t \in \mathbb{T}$, in general
$$
\left ( \int_{a}^{t}h(\tau) \diamondsuit_{\alpha} \tau
\right)^{\diamondsuit_{\alpha}} \ne h(t) \, .
$$
Although the fundamental theorem of calculus
does not hold for the $\diamondsuit_{\alpha}$-integral,
other properties do. Let $a, b, t \in \mathbb{T}$, $c \in \mathbb{R}$. Then,
\begin{itemize}
\item[(i)] $\int_{a}^{t}\{ f(\tau)+g(\tau) \} \diamondsuit_{\alpha} \tau =
\int_{a}^{t} f(\tau) \diamondsuit_{\alpha} \tau + \int_{a}^{t}
g(\tau) \diamondsuit_{\alpha} \tau$;

\item[(ii)] $\int_{a}^{t} c f(\tau) \diamondsuit_{\alpha} \tau =
c \int_{a}^{t} f(\tau) \diamondsuit_{\alpha} \tau$;

\item[(iii)] $\int_{a}^{t} f(\tau) \diamondsuit_{\alpha} \tau =
\int_{a}^{b} f(\tau) \diamondsuit_{\alpha} \tau + \int_{b}^{t}
f(\tau) \diamondsuit_{\alpha} \tau$;

\item[(iv)] If $f(t)\geq 0$ for all $t$, then $\int_a^b
f(t)\Diamond_\alpha t\geq 0$;

\item[(v)] If $f(t)\leq g(t)$ for all $t$, then $\int_a^b
f(t)\Diamond_\alpha t\leq\int_a^b g(t)\Diamond_\alpha t$;

\item[(vi)] If $f(t)\geq 0$ for all $t$, then 
$f(t) \equiv 0$ if and only if $\int_a^b f(t)\Diamond_\alpha t=0$.
\end{itemize}

In \cite{srd} a diamond-$\alpha$ Jensen's inequality 
on time scales is proved: let $c$ and $d$ be real numbers,
$u:[a,b]\rightarrow(c,d)$ be continuous and
$h:(c,d)\rightarrow\mathbb{R}$ be a convex function. Then,
\begin{equation}
\label{eq:JI}
h \left( \frac{\int_{a}^{b}u(t) \diamondsuit_{\alpha} t}{b-a}\right
) \leq \frac{1}{b-a} \int_{a}^{b}h(u(t)) \diamondsuit_{\alpha} t \, .
\end{equation}
We will use the diamond-$\alpha$ Jensen's inequality 
\eqref{eq:JI} in the proof of our Theorem~\ref{thm33}.

%---------------------------------------------------

\section{Main Results}
\label{sec:MR}

In the sequel we assume that $\mathbb{T}$ is a time scale,
$a$, $b \in \mathbb{T}$ with $a < b$, and $[a,b]$ denote
$[a,b] \cap \mathbb{T}$. We begin by proving some auxiliary lemmas.

\begin{lemma}\label{lem31}
If $f \in  \mathbb{C}([a, b], \mathbb{R})$ 
satisfy $m \leq f(t) \leq M$ for all $t \in [a,b]$ and
$
\int_{a}^{b}f(t)\diamondsuit_{\alpha}t =0 \, ,
$
then
$$
\frac{1}{b-a}\int_{a}^{b}f^{2}(t)\diamondsuit_{\alpha}t  \leq
\frac{1}{4}(M-m)^{2}.
$$
\end{lemma}

\begin{proof}
Define $\varphi(t)= \frac{f(t)-m}{M-m}$. It follows that
$f(t)= m +(M-m)\varphi(t)$. One can easily check that $0 \leq
\varphi(t) \leq 1$. Then,
\begin{equation*}
\begin{split}
\frac{1}{b-a} &\int_{a}^{b}f^{2}(t)\diamondsuit_{\alpha}t \\
& = \frac{1}{b-a} \int_{a}^{b} \left( m+(M-m)\varphi(t)
\right)^{2}\diamondsuit_{\alpha}t \\
& =\frac{1}{b-a}\int_{a}^{b}  \left(m^{2}+2m(M-m)\varphi(t)+
(M-m)^{2}\varphi^{2}(t)\right)\diamondsuit_{\alpha}t\\
& = m^{2}+ \frac{2m(M-m)}{b-a} \int_{a}^{b} \varphi(t)
\diamondsuit_{\alpha}t
+ \frac{(M-m)^{2}}{b-a} \int_{a}^{b} \varphi^{2}(t) \diamondsuit_{\alpha}t\\
& \leq m^{2}+ \frac{2m(M-m)}{b-a} \int_{a}^{b} \varphi(t)
\diamondsuit_{\alpha}t
+ \frac{(M-m)^{2}}{b-a} \int_{a}^{b} \varphi(t) \diamondsuit_{\alpha}t\\
&= m^{2} + \frac{2m(M-m)}{b-a} \int_{a}^{b} \frac{f(t)-m}{M-m}
\diamondsuit_{\alpha}t+ \frac{(M-m)^{2}}{b-a} \int_{a}^{b}\frac{f(t)-m}{M-m}  \diamondsuit_{\alpha}t\\
&=m^{2} + 2m(M-m)\left(\frac{-m}{M-m}\right)+ (M-m)^{2}\left(\frac{-m}{M-m}\right)\\
& = -mM = \frac{1}{4}\left( (M-m)^{2}-(M+m)^{2}\right)\\
& \leq \frac{1}{4} (M-m)^{2}.
\end{split}
\end{equation*}
\end{proof}

\begin{lemma}\label{lem32}
If $f \in  \mathbb{C}([a, b], \mathbb{R})$ satisfy 
$m \leq f(t) \leq M$ for all $t \in [a,b]$ 
and $\int_{a}^{b}f(t)\diamondsuit_{\alpha}t \neq 0$, then
$$
\frac{1}{b-a}\int_{a}^{b}f^{2}(t)\diamondsuit_{\alpha}t-
\left ( \frac{1}{b-a}\int_{a}^{b}f(t)\diamondsuit_{\alpha}t
\right)^{2}
\leq \frac{1}{4}(M-m)^{2} \, .
$$
\end{lemma}

\begin{proof}
Setting
$$
\frac{1}{b-a}\int_{a}^{b}f(t)\diamondsuit_{\alpha}t= I(b-a) \, ,
$$
$I \in \mathbb{R}$, and $F(t)= f(t)-I(b-a)$, we then have
$$
m- I(b-a) \leq F(t) \leq  M- I(b-a).
$$
Therefore,
\begin{equation}\label{eqI}
\begin{split}
\frac{1}{b-a}\int_{a}^{b}F(t)\diamondsuit_{\alpha}t &=
\frac{1}{b-a}\int_{a}^{b}(f(t)-I(b-a))\diamondsuit_{\alpha}t\\
& = \frac{1}{b-a}\int_{a}^{b}f(t)\diamondsuit_{\alpha}t -I(b-a)\\
& = 0.
\end{split}
\end{equation}

It follows from Lemma~\ref{lem31} that
\begin{equation}\label{eq1}
\begin{split}
\frac{1}{b-a}\int_{a}^{b}F^{2}(t)\diamondsuit_{\alpha}t  &
\leq \frac{1}{4} \left( (M-I(b-a))-(m-I(b-a)) \right)^{2}\\
 & = \frac{1}{4}(M-m)^{2}.
\end{split}
\end{equation}

On the other hand, using \eqref{eqI} we have
\begin{equation}\label{eq2}
\begin{split}
\frac{1}{b-a} &\int_{a}^{b}F^{2}(t)\diamondsuit_{\alpha}t- 
\left(\frac{1}{b-a}\int_{a}^{b}F(t)\diamondsuit_{\alpha}t \right)^{2} \\
& = \frac{1}{b-a}\int_{a}^{b}F^{2}(t)\diamondsuit_{\alpha}t\\
& =\frac{1}{b-a}\int_{a}^{b} \left (f^{2}(t)-2 I(b-a)f(t)+
I^{2}(b-a) \right ) \diamondsuit_{\alpha}t\\
& = \frac{1}{b-a}\int_{a}^{b} f^{2}(t) \diamondsuit_{\alpha}t -2
I^{2}(b-a)+I^{2}(b-a)\\
& = \frac{1}{b-a}\int_{a}^{b} f^{2}(t) \diamondsuit_{\alpha}t - I^{2}(b-a)\\
& = \frac{1}{b-a}\int_{a}^{b} f^{2}(t) \diamondsuit_{\alpha}t -
\left ( \frac{1}{b-a}\int_{a}^{b}f(t)\diamondsuit_{\alpha}t
\right)^{2}.
\end{split}
\end{equation}

Then, using \eqref{eq1} and \eqref{eq2}, we conclude with
$$
\frac{1}{b-a}\int_{a}^{b} f^{2}(t) \diamondsuit_{\alpha}t - \left (
\frac{1}{b-a}\int_{a}^{b}f(t)\diamondsuit_{\alpha}t \right)^{2} \leq
\frac{1}{4}(M-m)^{2}.
$$
\end{proof}

From Lemma~\ref{lem31} and Lemma~\ref{lem32} we have the
following corollary.

\begin{corollary}\label{cor33}
If $f \in  \mathbb{C}([a, b], \mathbb{R})$ 
satisfy $m \leq f(t) \leq M$ for all $t \in [a,b]$, then
$$
\frac{1}{b-a}\int_{a}^{b}f^{2}(t)\diamondsuit_{\alpha}t-
\left ( \frac{1}{b-a}\int_{a}^{b}f(t)\diamondsuit_{\alpha}t
\right)^{2} \leq \frac{1}{4}(M-m)^{2}.
$$
\end{corollary}

We are now in conditions to prove the intended
Gr\"{u}ss inequality.

\begin{theorem}\label{thm33} (The diamond-$\alpha$ Gr\"{u}ss 
inequality on time scales, $\alpha \in [0,1]$).
Let $\mathbb{T}$ be a time scale, $a, b \in \mathbb{T}$ with $a < b$. 
If $f, g \in \mathbb{C}([a, b], \mathbb{R})$ satisfy
$\varphi \leq f(t) \leq \Phi$ and $\gamma \leq g(t) \leq \Gamma$ 
for all $t \in [a,b] \cap \mathbb{T}$, then
\begin{equation}
\label{eq3}
\left|\frac{1}{b-a}\int_{a}^{b}f(t) g(t)\diamondsuit_{\alpha}t -
\frac{1}{(b-a)^{2}}\int_{a}^{b}f(t) \diamondsuit_{\alpha}t
 \int_{a}^{b}g(t) \diamondsuit_{\alpha}t  \right|
  \leq
\frac{1}{4}(\Phi-\varphi)(\Gamma-\gamma).
\end{equation}
\end{theorem}

\begin{proof}
A straightforward computation leads to
\begin{equation}\label{equation5}
\begin{split}
 &\frac{1}{b-a}\int_{a}^{b}f(t) g(t)\diamondsuit_{\alpha}t -
\frac{1}{(b-a)^{2}}\int_{a}^{b}f(t) \diamondsuit_{\alpha}t
 \int_{a}^{b}g(t) \diamondsuit_{\alpha}t \\
 & = \frac{1}{4} \left ( \frac{1}{(b-a)}\int_{a}^{b}\left ((f(t)+g(t))^{2}-(f(t)-g(t))^{2}
  \right ) \diamondsuit_{\alpha}t \right. \\
  & \qquad - \left. \frac{4}{(b-a)^{2}} \int_{a}^{b}f(t) \diamondsuit_{\alpha}t
 \int_{a}^{b}g(t) \diamondsuit_{\alpha}t    \right ).
\end{split}
\end{equation}
If we consider the function $h(x)=x^{2}$, which is
obviously convex, then using the diamond-$\alpha$ 
Jensen's inequality on time scales \eqref{eq:JI}, we obtain
\begin{equation}\label{eq4}
\left( \frac{\int_{a}^{b} (f(t)-g(t))
\diamondsuit_{\alpha}t}{b-a}\right )^{2} \leq  \frac{1}{b-a}
\int_{a}^{b} (f(t)-g(t))^{2}\diamondsuit_{\alpha} t.
\end{equation}
 Then we have by \eqref{equation5} and \eqref{eq4} that
\begin{equation}\label{equation7}
\begin{split}
 \frac{1}{b-a} & \int_{a}^{b}f(t) g(t)\diamondsuit_{\alpha}t -
\frac{1}{(b-a)^{2}}\int_{a}^{b}f(t) \diamondsuit_{\alpha}t
 \int_{a}^{b}g(t) \diamondsuit_{\alpha}t \\
 &\leq  \frac{1}{4}\Biggl\{ \frac{1}{(b-a)}\int_{a}^{b}\left
((f(t)+g(t)\right)^{2}\diamondsuit_{\alpha} t -
\frac{1}{(b-a)^{2}}\left( \int_{a}^{b} (f(t)-g(t))
\diamondsuit_{\alpha}t \right)^{2} \\
& \qquad -\frac{4}{(b-a)^{2}} \int_{a}^{b}f(t)
\diamondsuit_{\alpha}t
 \int_{a}^{b}g(t) \diamondsuit_{\alpha}t \Biggr\}\\
 &\leq  \frac{1}{4}\Biggl\{ \frac{1}{(b-a)}\int_{a}^{b}\left
((f(t)+g(t)\right)^{2}\diamondsuit_{\alpha} t- \frac{1}{(b-a)^{2}}
\left( \int_{a}^{b} (f(t)+g(t)) \diamondsuit_{\alpha}t \right)^{2} \\
& \qquad  +\frac{1}{(b-a)^{2}} \left( \int_{a}^{b} (f(t)+g(t)) 
\diamondsuit_{\alpha}t \right)^{2}- \frac{1}{(b-a)^{2}}\left(
\int_{a}^{b} (f(t)-g(t)) \diamondsuit_{\alpha}t \right)^{2} \\
& \qquad   -\frac{4}{(b-a)^{2}} \int_{a}^{b}f(t)
\diamondsuit_{\alpha}t
 \int_{a}^{b}g(t) \diamondsuit_{\alpha}t \Biggr\} \\
& \leq  \frac{1}{4}\left\{ \frac{1}{(b-a)}\int_{a}^{b}\left
((f(t)+g(t)\right)^{2}\diamondsuit_{\alpha} t- \frac{1}{(b-a)^{2}}
\left( \int_{a}^{b} (f(t)+g(t)) \diamondsuit_{\alpha}t \right)^{2}
\right \},
\end{split}
\end{equation}
since
\begin{multline*}
\frac{1}{(b-a)^{2}} \left( \int_{a}^{b} (f(t)+g(t))
\diamondsuit_{\alpha}t \right)^{2}- \frac{1}{(b-a)^{2}}\left(
\int_{a}^{b} (f(t)-g(t)) \diamondsuit_{\alpha}t \right)^{2}\\
 -\frac{4}{(b-a)^{2}} \int_{a}^{b}f(t) \diamondsuit_{\alpha}t
 \int_{a}^{b}g(t) \diamondsuit_{\alpha}t =0.
\end{multline*}

On the other hand, we have $\varphi+\Phi \leq (f+g)(t) \leq
\gamma+\Gamma$.  Applying Corollary~\ref{cor33} to the function
$f+g$, we get
\begin{multline}
\label{eq5}
\frac{1}{b-a}\int_{a}^{b}(f(t)+g(t))^{2} \diamondsuit_{\alpha} t 
- \frac{1}{(b-a)^{2}} \left(\int_{a}^{b}(f(t)+g(t))
\diamondsuit_{\alpha} t \right)^{2} \diamondsuit_{\alpha} t \\
\leq \frac{1}{4}((\Phi+\Gamma) - (\varphi+\gamma))^{2}=
\frac{1}{4}(\Phi+\Gamma-\varphi-\gamma)^{2}.
\end{multline}
Gathering \eqref{eq4}, \eqref{equation7} and \eqref{eq5}, we obtain:
\begin{equation*}
\begin{split}
\frac{1}{b-a} & \int_{a}^{b}f(t) g(t)\diamondsuit_{\alpha}t -
\frac{1}{(b-a)^{2}}\int_{a}^{b}f(t) \diamondsuit_{\alpha}t
 \int_{a}^{b}f(t) \diamondsuit_{\alpha}t  \\
& \leq \frac{1}{16} \left((\Phi+\Gamma)-(\varphi+\gamma)\right)^{2}\\
& \leq \frac{1}{4} \left(
(\Phi-\varphi)(\Gamma-\gamma)+\frac{1}{4}(\Phi-\varphi-(\Gamma-\gamma))^{2}
\right).
\end{split}
\end{equation*}
We consider now two cases: (i) if
$\Phi-\varphi= \Gamma-\gamma$, then
\begin{equation}
\label{eq6}
 \frac{1}{b-a}\int_{a}^{b}f(t)
g(t)\diamondsuit_{\alpha}t - \frac{1}{(b-a)^{2}}\int_{a}^{b}f(t)
\diamondsuit_{\alpha}t
 \int_{a}^{b}g(t) \diamondsuit_{\alpha}t
 \leq \frac{1}{4} (\Phi-\varphi)(\Gamma-\gamma) \, ;
\end{equation}
(ii) if $\Phi-\varphi\neq  \Gamma-\gamma$, let us define
 $$
 \beta = \sqrt{\frac{\Gamma-\gamma}{\Phi-\varphi}}\, , \quad
\mu =  \sqrt{\frac{\Phi-\varphi}{\Gamma-\gamma}},
$$
and
$$
f_{1}(t) = \beta f(t) \, , \quad g_{1}(t)= \mu g(t).
$$
Note that $\beta \mu =1$. Then, we have
$$
\overline{m}_{1}= \beta \varphi \leq f_{1}(t) \leq \beta \Phi=
\overline{M}_{1}, \, \overline{m}_{2}= \mu \gamma \leq g_{1}(t) \leq
\mu \Gamma= \overline{M}_{2} ,
$$
and it follows that
\begin{equation}\label{eqii}
\begin{split}
\overline{M}_{1}-\overline{m}_{1}&= \beta \Phi -\beta \varphi= \beta
( \Phi-\varphi)\\
&= \sqrt{(\Phi-\varphi)(\Gamma-\gamma)} = \mu(\Gamma-\gamma)=
\overline{M}_{2}-\overline{m}_{2}.
\end{split}
\end{equation}
Using the fact that $\beta \mu =1$, \eqref{eq6} and
\eqref{eqii} (with $f_{1}$, $g_{1}$), we get
\begin{multline}
\label{eq9}
\frac{1}{b-a}\int_{a}^{b}  f(t)
g(t)\diamondsuit_{\alpha}t - \frac{1}{(b-a)^{2}}\int_{a}^{b} f(t)
\diamondsuit_{\alpha}t
 \int_{a}^{b}  g(t) \diamondsuit_{\alpha}t \\
= \frac{1}{b-a}\int_{a}^{b} \beta \mu f(t)
g(t)\diamondsuit_{\alpha}t - \frac{1}{(b-a)^{2}}\int_{a}^{b} \beta
f(t) \diamondsuit_{\alpha}t
 \int_{a}^{b} \mu g(t) \diamondsuit_{\alpha}t \\
 = \frac{1}{b-a}\int_{a}^{b}f_{1}(t)
g_{1}(t)\diamondsuit_{\alpha}t -
\frac{1}{(b-a)^{2}}\int_{a}^{b}f_{1}(t) \diamondsuit_{\alpha}t
 \int_{a}^{b}g_{1}(t) \diamondsuit_{\alpha}t \\
 \leq \frac{1}{4} (\overline{M}_{1}-\overline{m}_{1})
(\overline{M}_{2}-\overline{m}_{2})  = \frac{1}{4} \beta \mu
(\Phi-\varphi)(\Gamma-\gamma)
= \frac{1}{4} (\Phi-\varphi)(\Gamma-\gamma).
\end{multline}
If we now consider the case of $-f$ we have $-\Phi \leq -f(t) \leq
-\varphi$. Using \eqref{eq9},
\begin{multline}
\label{eq10}
 \frac{1}{b-a}\int_{a}^{b}  (-f)(t)
g(t)\diamondsuit_{\alpha}t - \frac{1}{(b-a)^{2}}\int_{a}^{b} (-f)(t)
\diamondsuit_{\alpha}t
 \int_{a}^{b}  g(t) \diamondsuit_{\alpha}t \\
= -\left \{ \frac{1}{b-a}\int_{a}^{b}  f(t)
g(t)\diamondsuit_{\alpha}t - \frac{1}{(b-a)^{2}}\int_{a}^{b} f(t)
\diamondsuit_{\alpha}t
 \int_{a}^{b}  g(t) \diamondsuit_{\alpha}t \right \} \\
 \leq \frac{1}{4} (-\varphi- (-\Phi))(\Gamma-\gamma)
 = \frac{1}{4} (\Phi-\varphi)(\Gamma-\gamma).
\end{multline}

Then, using \eqref{eq9} and \eqref{eq10}, we arrive to
\begin{equation*}
\left| \frac{1}{b-a}\int_{a}^{b}  f(t) g(t)\diamondsuit_{\alpha}t -
\frac{1}{(b-a)^{2}}\int_{a}^{b} f(t) \diamondsuit_{\alpha}t
 \int_{a}^{b}  g(t) \diamondsuit_{\alpha}t \right|
\leq \frac{1}{4} (\Phi-\varphi)(\Gamma-\gamma).
\end{equation*}
This completes the proof of our Theorem~\ref{thm33}.
\end{proof}

\begin{remark}
When $\alpha =1$, we have the following   
Gr\"{u}ss inequality on time scales:
\begin{equation*}
\left| \frac{1}{b-a}\int_{a}^{b}  f(t) g(t)\Delta t -
\frac{1}{(b-a)^{2}}\int_{a}^{b} f(t) \Delta t
 \int_{a}^{b}  g(t) \Delta t \right|
\leq \frac{1}{4} (\Phi-\varphi)(\Gamma-\gamma).
\end{equation*}
\end{remark}

\begin{remark}
For $\mathbb{T}=\mathbb{R}$ Theorem~\ref{thm33}
gives the classical inequality \eqref{eq:class:Gruss}.
\end{remark}

Applying the Gr\"{u}ss inequality in Theorem~\ref{thm33} 
to the time scale $\mathbb{T}= \mathbb{Z}$ with
$a=0$, $b=n$, and $f(i)=x_{i}$, $i = 1,\ldots, n$,
we arrive to the following corollary which improves \cite{dragomir}.

\begin{corollary}
If
$\varphi \leq x_{i} \leq  \Phi$ and $\gamma \leq y_{i} \leq  \Gamma$ 
for all $1 \leq i \leq n$, then the following discrete  
Gr\"{u}ss inequality holds:
$$
\left| \frac{1}{n} \sum_{i=1}^{n}x_{i}y_{i} -\frac{1}{n^{2}}
\sum_{i=1}^{n}x_{i} \sum_{i=1}^{n}y_{i} \right| \leq \frac{1}{4}
(\Phi-\varphi)(\Gamma-\gamma) \, .
$$
\end{corollary}

Theorem~\ref{thm33} is valid for an arbitrary time scale.
For example, let $q > 1$ and $\mathbb{T}= \{ q^{N_{0}} \}$.

\begin{corollary}(Quantum Gr\"{u}ss inequality)
If $f$ and $g$ satisfy $\varphi \leq f(q^{i}) \leq  \Phi$ 
and $\gamma \leq g(q^{i}) \leq  \Gamma$ 
for all $q^i$, $i=m,\ldots,n$, then the following inequality holds:
 \begin{multline*}
\left | \frac{\sum_{i=m}^{n-1}
q^{i}f(q^{i+1})g(q^{i+1})}{\sum_{i=m}^{n-1}q^{i}} -\frac{1}{\left(
\sum_{i=m}^{n-1} q^{i} \right)^{2}} \left( \sum_{i=m}^{n-1}
q^{i}f(q^{i+1}) \right) \left( \sum_{i=m}^{n-1}
q^{i}g(q^{i+1}) \right) \right |\\
\leq \frac{1}{4} (\Phi-\varphi)(\Gamma-\gamma).
\end{multline*}
\end{corollary}

%---------------------------------------------------

%---------------------------------------------------

\end{document}